\documentclass[letterpaper, 10 pt, conference]{ieeeconf}

\IEEEoverridecommandlockouts 
\usepackage{stackrel}
\usepackage{graphicx}
\usepackage{siunitx}
\graphicspath{{./}}
\usepackage{amssymb,amsmath}
\usepackage{multicol}
\usepackage{bm} 

\newtheorem{definition}{Definition}
\newtheorem{lemma}{Lemma}
\newtheorem{theorem}{Theorem}
\newtheorem{corollary}{Corollary}

\newtheorem{proofoflemma}{Proof of Lemma}
\newtheorem{proofoftheorem}{Proof of Theorem}
\title{\LARGE \bf Heterogeneous Stochastic Interactions for Multiple Agents in a Multi-armed Bandit Problem}

\author{Udari Madhushani
and Naomi Ehrich Leonard% <-this % stops a space
	\thanks{This research has been supported in part by ARO grant W911NF-18-1-0325 and ONR grant N00014-14-1-0635. Department of Mechanical and Aerospace Engineering, Princeton University, NJ 08544, USA.
		{\tt\small \{udarim,naomi\}@princeton.edu}}%
}

\begin{document}

\maketitle
\thispagestyle{empty}
\pagestyle{empty}

\begin{abstract}
We define and analyze a multi-agent multi-armed bandit problem in which decision-making agents can observe the choices and rewards of their neighbors.  Neighbors are defined by a network graph with heterogeneous and stochastic interconnections. These interactions are determined by the “sociability” of each agent, which corresponds to the probability that the agent observes its neighbors.  We design an algorithm for each agent to maximize its own expected cumulative reward and prove performance bounds that depend on the sociability of the agents and the network structure.  We use the bounds to predict the rank ordering of agents according to their performance and verify the accuracy analytically and computationally.
\end{abstract}
\allowdisplaybreaks

\section{Introduction} \label{sect:Introduction}
%\cite{Baltaoglu,Lai}

Animal and robotic foraging problems that involve searching over spatially distributed patches with uncertain distribution of resource (reward) result in the explore-exploit dilemma. Each animal or robotic forager chooses which patch to sample, sequentially in time. The dilemma at every time in the sequence is whether to choose a well sampled patch that is expected to reap the highest reward (exploit) or to choose a poorly sampled patch that is expected to reduce uncertainty and possibly identify an option with an even higher reward (explore). Choices to explore can be costly or risky since they may not return much reward, but successful exploiting requires sufficient information that comes from exploration. The challenge is sorting out protocols for when and where to explore versus when and where to exploit. 

In decision theory, multi-armed bandit (MAB) problems serve as models that capture the salient features of the explore-exploit trade-off \cite{Sutton,Robbins}. Thus, advances in addressing MAB problems directly benefit understanding and solving foraging and search problems. The MAB problem is analogous to a scenario in which an agent is repeatedly faced with several different options, each returning uncertain reward, and aims to make a sequence of choices to maximize the cumulative reward \cite{Gittins}. This is equivalent to minimizing the cumulative regret \cite{LaiRobbins}. 

% Over decades, optimal strategies have been developed to realize the MAB objective. In the standard MAB problem the reward distributions are stationary. Thus, if the mean values of the reward distributions of all the options are known to the agent, then in order to maximize the cumulative reward, the agent only has to sample the option with the maximum mean. In reality the means are not available to the agent.  So the agent should balance choosing options that will provide information to better estimate the mean values with choosing options that it expects will maximize the cumulative reward. This is called the explore-exploite dilemma.  

In their seminal work, Lai and Robbins \cite{LaiRobbins} established a lower bound for the expected cumulative regret in the finite time horizon case. Specifically, they derived a logarithmic lower bound for the expected number of times a sub-optimal option needs to be sampled by an optimal sampling rule. They also established a confidence bound and a sampling rule to achieve logarithmic cumulative regret. These results were further simplified in \cite{AgrawalSimpl} by introducing a confidence bound using a sample mean based method.  
Improving on these results, a family of Upper Confidence Bound (UCB) algorithms for achieving asymptotic and uniform logarithmic cumulative regret were proposed in \cite{Auer}. These UCB algorithms are based on the notion that %the desired goal of 
minimizing the expected cumulative regret is realized by choosing an appropriate uncertainty model, which results in an optimal trade-off between reward gain and information gain through uncertainty. 

For the standard MAB problem the reward associated with an option is considered  to be an iid stochastic process. Therefore, in the frequentist setting, the natural way of estimating the expectation of the reward is to consider the sample average \cite{LaiRobbins,AgrawalSimpl,Auer}. The papers \cite{Kauffman,Reverdy} present how to incorporate prior knowledge about reward expectation in the estimation step by leveraging the theory of conditional expectation in the Bayesian setting.

The papers \cite{1104491,PeterCDC} extend this work to the multi-agent setting where the explore-exploit problem is defined for a group of agents. The objective is to understand how the performance of individuals or the group can benefit from inter-agent observations. A centralized multi-agent setting is considered in  \cite{1104491} and a decentralized setting is considered in \cite{6763073}. The papers \cite{PeterECC,PeterCDC} use a running consensus algorithm in which agents observe the estimates of their neighbors. A multi-agent multi-armed bandit (MAMAB) problem, where agents observe instantaneous rewards and actions in a leader-follower setting, is considered in \cite{Gopalan16,LandgrenCDC2018}.

In the above works, agents observe their neighbors, defined according to a static network graph. In this paper we consider a MAMAB problem where agents observe instantaneous rewards and actions of their neighbors through stochastic interactions. Observing rewards and actions rather than estimates is motivated by the foraging success of social groups where agents can observe neighbors even when their neighbors don't necessarily want to share what they know \cite{Tilman2016MaintainingCI}. % but may not be inclined to broadcast what they know.
%explicitly communicate with their neighbours, but they observe their neighbour. 
%Observing does not rely on cooperation between neighbours.  
Further, we assume that each agent $k$ only observes its neighbors with some probability $p_k$.   
This provides a framework for evaluating efficiency and robustness to changes in communication in terms of the $p_k$.
%e.g., some agents with lower $p_k$, as well as robustness to loss of reliability in communication, e.g., a drop in some $p_k$. %For instance a group of animals who observe each other without explicitly signalling can demonstrate good performance in social foraging \cite{}. 
%Another example can be given by considering the performance of a group of fishermen who do not broadcast but observe their neighbours \cite{}. 

The setting is formulated by defining an underlying undirected network graph and imposing directed observation probabilities on each edge. The underlying graph models the inherent observation constraints present in the network. Imposed observation probabilities capture the heterogeneous social effort of agents in observing neighbors. We introduce the notion of \textit{sociability} of each agent $k$ as the likelihood $p_k$ that the agent observes its neighbors. We derive analytical upper bounds for expected cumulative regret and propose a measure to rank agents according to their relative performance as a function of each agent's sociability and the sociability of its neighbors. We show that our model predicts how high performance requires an agent to have {\em both} high sociability {\em and} neighbors with low sociability. This is an important result: it implies that making an investment to observe neighbors may be worthwhile only if those neighbors are sufficiently explorative.
% of an agent's neighbors on the performance of the agent.

%The rest of the paper is organized as follows. 
In Section \ref{Secn:MAMAB} we introduce the MAMAB problem with time-varying (stochastic) observation structure. We propose an efficient sampling rule for an agent to maximize cumulative expected reward. We analyze the performance of the proposed sampling rule in Section \ref{Secn:Regret}. In Section \ref{SubSecn:Regret} we analytically upper bound the expected cumulative regret and in Section \ref{SubSecn:Perfom} we propose a measure to predict the ranks of agents according to their relative performance. In Section \ref{Secn:Simu} we provide numerical simulation results and computationally validate the analytical results. %Simulation results for an all-to-all observation network graph and a cyclic observation network graph are given in Section \ref{SubSecn:AlltoAll} and Section \ref{SubSecn:Cyclic} respectively. 
We conclude in Section \ref{Secn:Concl} and provide additional mathematical details in the Appendix.

%%%%%%%%%%%%%%%%%%%%%%%%%%%%%%%%%%%%%%%%%%%

\section{Multi-agent Multi-armed Bandit Problem}\label{Secn:MAMAB}

Consider a MAMAB with $K$ agents and $N$ arms, which represent $N$ options (patches).  Let $X_{i}$ be a sub-Gaussian random variable with variance proxy $(\sigma_i^{\prime})^2$ that denotes the reward associated with option $i\in \mathcal{I} \triangleq \{1,2,\ldots,N\}.$ Define $\mu_{i}$ as the expected reward $\mathbb{E}(X_{i})$ of option $i$. Let $i^{*}$ be the optimal option such that $\mu_{i^{*}}=\max \{\mu_{1},\ldots,\mu_{N}\}$.
Each agent $k \in \{1,\ldots, K\}$ chooses one option at each time step with the goal of maximizing its cumulative reward, which is equivalent to minimizing its cumulative regret. We assume that the $\mu_{i}$  are unknown and the $\sigma_{i}^{\prime}$ are known to the agents  (e.g., if $\sigma_{i}^{\prime}$ is the magnitude of sensor noise). 

Let $\mathcal{G}(\mathcal{V},\mathcal{E})$  be an undirected graph that encodes the observation structure of the system. If there exists an observation link $e_{kj}$ from agent $k$ to agent $j$, $j,k \in \mathcal{V}$, % the $k^\mathrm{th}$ agent to the $j^\mathrm{th}$ agent, 
then  $\{k,j\},\{j,k\} \in \mathcal{E}$,  $e_{kj}= e_{jk}=1,$  and we say 
agents $k$ and $j$ are neighbors.
For social or robotic groups searching over physical space, observation links may exist between pairs of agents when they are sufficiently close to be visible (or otherwise observable) to each other.
%Let $e_{kj}$ be the observation link from the $k^\mathrm{th}$ agent to the $j^\mathrm{th}$ agent. If there exists an observation link between $k^{\mathrm{th}}$ agent and $j^{\mathrm{th}}$ agent then  $\{k,j\}\in \mathcal{E}$ and $e_{kj}= e_{jk}=1.$ 
%In the examples of animal and fishermen groups this corresponds to the range of visibility. Individuals can observe others who are within their range of visibility. 
%In this case we define agents $k$ and $j$ to be neighbors. 
Let $d_{k}$ be the number of neighbors  and $\mathcal{V}_{d_{k}}$  the set of neighbors of agent $k$. 

Let $\varphi_{k}^{t}\in \mathcal{I}$ and $X_i^t$ be random variables that denote the option chosen and reward received by agent $k$ at time $t$, respectively. Let $X_i^t$, $\forall t$, be i.i.d copies of $X_i.$
%the reward received by  agent $k$ at time $t$. 
Define $\sigma_i=\sqrt{K}\sigma_i^{\prime}.$
Consider the probability space $(\Omega,\mathcal{U},\mathbb{P})$ and the increasing sequence of subalgebras $\mathcal{F}_{0}\subset\mathcal{F}_{1}\cdots \subset\mathcal{F}_{t}\cdots \subset\mathcal{F}_{T-1}\subset \mathcal{U}$ for $t=0,1,\cdots,T$, where 
$\mathbb{P}$ is the probability measure on the sigma algebra $\mathcal{U}$ of $\Omega$. Here $\mathcal{F}_{t}$ is the sigma algebra generated by information available at time $t$. Let $\mathbb{I}_{\{\varphi_{k}^{t} =i\}}$ be a $\mathcal{F}_{t-1}$ measurable indicator random variable that takes the value one if  option $i$ is chosen by agent $k$ at time $t$ and is zero otherwise. Define $\mathbb{I}^{t}_{\{k,j\}}$ to be an $\mathcal{F}_{t-1}$ measurable indicator  random variable that takes the value one if agent $k$ can observe agent $j$ at time $t$ and is zero otherwise. 

Let $p_k$ be the probability that agent $k$ observes the instantaneous actions and rewards of its neighbors. Then we have $\mathbb{E}\left(\mathbb{I}^{t}_{\{k,j\}}\right)=p_{k},\forall j\neq k$ such that $\{k,j\}\in\mathcal{E}.$  We let $\mathbb{I}^{t}_{\{k,k\}}=1$ and $\mathbb{I}^{t}_{\{k,j\}}=0,\forall j\neq k$ such that $\{k,j\}\notin\mathcal{E}.$ An agent $k$ that has high probability $p_k$ is more likely to obtain observations from its neighbors. We introduce the notion of an agent's \textit{sociability} to refer to its value of $p_k$: we interpret agents with high $p_{k}$ values to be more sociable and agents with low $p_{k}$ values to be less sociable.

In order to maximize the cumulative reward in the long run, agents need to both identify the best options through exploring and sample the best options through exploiting. Since agents with high sociability values are more likely to obtain a greater number of observations, they can identify the best options with less exploring. However the usefulness of their observations is affected by the sociability values of their neighbors. Better performance can be obtained by an agent when it observes neighbors that do a lot of exploring (because they are less sociable), as compared to when it observes agents that do a lot of exploiting (because they are more sociable). This is, the agent will be able to exploit more, without compromising performance, when it has neighbors that are less sociable.

Let the number of times agent $k$ samples option $i$ in $T$ trials be given by the $\mathcal{F}_{T-1}$ measurable random variable  $n_{i}^{k}(T)\triangleq \sum_{t=1}^T\mathbb{I}_{\{\varphi_{k}^{t} =i\}}$. And let the total number of times that  agent $k$ observes rewards from option $i$ be the $\mathcal{F}_{T-1}$ measurable random variable $N_{i}^{k}(T)$, given as
\begin{align*}
N_{i}^{k}(T)=\sum_{t=1}^{T}\sum_{j=1}^{K}\mathbb{I}_{\{\varphi_{j}^{t} =i\}}\mathbb{I}^{t}_{\{k,j\}}.
\end{align*}
%which is the total number of times that  agent $k$ observes rewards from option $i$. 
We define a sampling rule for agent $k$ as follows. Let the $\mathcal{F}_{t}$ measurable random variable ${\hat{\mu}}_{i}^{k}(T)$ be the estimate of option $i$ by agent $k$ at time $t$. We define
\begin{align*}
\hat{\mu}_{i}^{k}(T)=\frac{S_{i}^{k}(T)}{N_{i}^{k}(T)},
\end{align*}
where $S_{i}^{k}(T)=\sum_{t=1}^{T}\sum_{j=1}^{K}X_{i}^t\mathbb{I}_{\{\varphi_{j}^{t} =i\}}\mathbb{I}^{t}_{\{k,j\}}$ is the total reward observed by agent $k$ from option $i$ in $T$ trials.  
\begin{definition}\label{def:UCB based}
%Let $\hat{\mu}_i^k(t)$ be an $\mathcal{F}_{t}$ measurable random variable which models the estimator of the $k^{\mathrm{th}}$ agent of $i^{\mathrm{th}}$ optation at $t^{\mathrm{th}}$ time step. 
The {\em sampling rule} $\{\varphi_{k}^{t}\}_1^{T}$ {\em for agent $k$ at time} $t \in \{1, \ldots, T\}$ is defined as
	\begin{align}
	\mathbb{I}_{\{\varphi_{k}^{t+1}=i\}}=\left\{
	\begin{array}{cl} 1 &, \:\:\:Q^{k}_i(t)=\max\{Q^{k}_1(t),\cdots,Q^{k}_N(t)\}\label{eq:UCBallocation}\\
  	0 &, \:\:\: {\mathrm{o.w.}}\end{array}\right.
	\end{align}
	with 	
	\begin{align}
	Q^{k}_i(t)&\triangleq \hat{\mu}_{i}^{k}(t)+C_i^k(t)\label{eq:UCBQ}\\
	C_i^k(t)&\triangleq\sigma_{i}\sqrt{\frac{2(\xi+1)\left(N_{i}^{k}(t)+f(t)\right)}{N_{i}^{k}(t)}\frac{\log t}{N_{i}^{k}(t)}}\label{eq:Uncertainity}
	\end{align}
	where $\xi>1$
and $f(t)$ a sublogarithmic nondecreasing nonnegative function.
\end{definition}

%%%%%%%%%%%%%%%%%%%%%%%%%%%%%%%%%%%%%%%
\section{Performance Analysis}\label{Secn:Regret}
In this section we proceed to analyze the performance of the proposed sampling rule by analyzing the expected cumulative regret of the agents.  Recall that the expected cumulative regret depends on the expected number of times the agents sample suboptimal options and the goal of each agent is to maximize its individual cumulative reward.

\subsection{Regret Analysis}\label{SubSecn:Regret}
Let $i$ be a suboptimal option. The total number of times agent $k$ samples from  option $i$ can be upper bounded as follows: 
\[
n_{i}^{k}(T)=\sum_{t=1}^T\mathbb{I}_{\{\varphi_{k}^{t} =i\}}
\leq \sum_{t=1}^{T}\mathbb{I}_{\{Q_{i}^{k}(t)\geq Q_{i^{*}}^{k}(t)\}}.
\]
%\begin{align*}
%n_{i}^{k}(T)&=\sum_{t=1}^T\mathbb{I}_{\{\varphi_{k}^{t} =i\}}\nonumber \\
%&\leq \sum_{t=1}^{T}\mathbb{I}_{\{Q_{i}^{k}(t)\geq Q_{i^{*}}^{k}(t)\}}
%\end{align*}
Here $\mathbb{I}_{\{Q_{i}^{k}(t)>Q_{i^{*}}^{k}(t)\}}$ is an indicator function that takes value one if the objective function of option $i$ is greater than the objective function of the optimal option $i^*$ and zero otherwise, i.e.,
\begin{align*}
\mathbb{I}_{\{Q_{i}^{k}(t)>Q_{i^{*}}^{k}(t)\}}=\left\{
	\begin{array}{cl} 1 &, \:\:\:Q^{k}_i(t)\geq Q^{k}_{i^{*}}(t)\\
  	0 &, \:\:\: {\mathrm{o.w.}}\end{array}\right.
\end{align*}
Thus we have 
\begin{align*}
\mathbb{E}\left(n_{i}^{k}(T)\right)\leq \sum_{t=1}^{T}\mathbb{P}\left(Q_{i}^{k}(t)\geq Q_{i^{*}}^{k}(t)\right). 
\end{align*}

Let $R_{i}^{k}(T)$ be the cumulative regret of agent $k$ from option $i$. Define $\Delta_{i}\triangleq \mu_{i^{*}}-\mu_{i}.$ Then we have, from \cite{Lai},
\begin{align}
\mathbb{E}\left(R_{i}^{k}(T)\right)=\Delta_{i}\mathbb{E}\left(n_{i}^{k}(T)\right).\label{eq:regretValue}
\end{align}
We proceed to analyze the expected number of samples from suboptimal options as follows. 

First we note that $\forall i,k,t$ we have 
\begin{align}
&\left\{Q_{i}^{k}(t)\geq Q_{i^{*}}^{k}(t)\right\}\subseteq \left\{\mu_{i^{*}}<\mu_{i}+2C_{i}^{k}(t)\right\}\nonumber\\
&\:\:\:\: \:\:\:\:\:\:\:\:\cup\left\{\hat{\mu}_{i^{*}}^{k}\leq \mu_{i^{*}}-C_{i^{*}}^{k}(t)\right\}\cup \left\{\hat{\mu}_{i}^{k}\geq \mu_{i}+C_{i}^{k}(t)\right\}\\
&\mathbb{E}\left(n_{i}^{k}(T)\right) \leq \sum_{t=1}^{T}\mathbb{P}\left(\mu_{i^{*}}<\mu_{i}+2C_{i}^{k}(t)\right)\nonumber\\
&+\sum_{t=1}^{T}\mathbb{P}\left(\hat{\mu}_{i^{*}}^{k}(t)\leq \mu_{i^{*}}-C_{i^{*}}^{k}(t)\right)+ \sum_{t=1}^{T}\mathbb{P}\left(\hat{\mu}_{i}^{k}(t)\geq \mu_{i}+C_{i}^{k}(t)\right)
.\label{eq:Prob}
\end{align}
% \begin{align}
% \left\{Q_{i}^{k}(t)\geq Q_{i^{*}}^{k}(t)\right\}\subseteq &\left\{\hat{\mu}_{i^{*}}^{k}\leq \mu_{i^{*}}-C_{i^{*}}^{k}(t)\right\}\cup \left\{\hat{\mu}_{i}^{k}\geq \mu_{i}-C_{i}^{k}(t)\right\}\nonumber\\
% &\:\:\:\cup\left\{\mu_{i^{*}}<\mu_{i}+2C_{i}^{k}(t)\right\}\nonumber\\
% \mathbb{E}\left(n_{i}^{k}(T)\right) \leq \sum_{t=1}^{T}\mathbb{P}&\left(\hat{\mu}_{i^{*}}^{k}\leq \mu_{i^{*}}-C_{i^{*}}^{k}(t)\right)+ \sum_{t=1}^{T}\mathbb{P}\left(\hat{\mu}_{i}^{k}\geq \mu_{i}-C_{i}^{k}(t)\right)\nonumber\\
% &\:\:\:+\sum_{t=1}^{T}\mathbb{P}\left(\mu_{i^{*}}<\mu_{i}+2C_{i}^{k}(t)\right).\label{eq:Prob}
% \end{align}

Next we proceed to analyze concentration probability bounds on the estimates of options. 
\begin{theorem}\label{thm:TailProb}
For any $\zeta>1$ and for $\sigma_i>0$ there exists a $\vartheta>0$ such that
\begin{align*}
    \mathbb{P}\left(\hat{\mu}_i^k(t)-{\mu}_i>\sqrt{\frac{\vartheta}{N_{i}^{k}(t)}}\right)\leq \frac{\nu\log Kt}{\exp(2\kappa \vartheta)}
\end{align*}
where 
\begin{align*}
\nu=\frac{1}{\log \zeta},\:\:\:\: \kappa=\frac{1}{\sigma_i^2\left(\zeta^{\frac{1}{4}}+\zeta^{-\frac{1}{4}}\right)^2}.
\end{align*}
\end{theorem}

{\em Proof of Theorem~\ref{thm:TailProb}}: See Appendix.

Using symmetry we conclude that
\begin{align*}
    \mathbb{P}\left(\Big |\hat{\mu}_i^k(t)-{\mu}_i\Big |>\sqrt{\frac{\vartheta}{N_{i}^{k}(t)}}\right)\leq \frac{\nu\log Kt}{\exp(2\kappa \vartheta)}.
\end{align*}

% Using theorem-\ref{thm:TailProb} we calculate explicit concentration probability bounds for the choice $C_i^k(t)$ given in equation (\ref{eq:Uncertainity})
\begin{lemma}\label{lem:TailProb}
For $\vartheta=2\sigma_{i}^{2}\left(\delta(\xi)+\delta'(\epsilon)\right)\log t$, with $\epsilon>0, \delta(\xi)>0, \delta'(\epsilon) = \delta(\xi)\epsilon/4$, there exists a $\zeta>1$ such that
\begin{align*}
    \mathbb{P}\left(\Big|\hat{\mu}_i^k(t)-{\mu}_i\Big |>\sigma_{i}\sqrt{\frac{2\left(\delta(\xi)+\delta'(\epsilon)\right)\log t}{N_{i}^{k}(t)}}\right)\leq \frac{\nu\log Kt}{t^{\delta(\xi)}}.
\end{align*}
\end{lemma}
\begin{proofoflemma}
Define $\Gamma(\zeta)=\frac{1}{\left(\zeta^{\frac{1}{4}}+\zeta^{-\frac{1}{4}}\right)^2}$ where $\zeta>1.$  $\Gamma(\zeta)$ is a monotonically decreasing function and
%\begin{align*}
$\lim_{\zeta\to 1}\Gamma(\zeta)=\frac{1}{4}$.
%\end{align*}
 This implies that $\kappa<\frac{1}{4\sigma_{i}^{2}}$ since $\kappa = \Gamma(\zeta)/\sigma_i^2$. Choose $\zeta > 1$ such that $\kappa=\frac{1}{(4+\epsilon)\sigma_{i}^{2}},$ where $\epsilon>0.$ 
 Then $\exists \delta'(\epsilon)>0$ such that $4\kappa(\delta(\xi)+\delta'(\epsilon))\sigma_{i}^{2}=\delta(\xi).$ 
 This proves that $\exists \zeta>1$ such that $2\kappa\vartheta=\delta(\xi)\log t$. The lemma follows from Theorem~\ref{thm:TailProb}. % it follows that
% \begin{align*}
 %   \mathbb{P}\left(\Big|\hat{\mu}_i^k(T)-{\mu}_i\Big |>\sigma_{i}\sqrt{\frac{2(\xi+1)\log T}{N_{i}^{k}(T)}}\right)\leq \frac{\nu\log KT}{T^{\xi+1}}.
%\end{align*}
%This concludes the proof of Lemma~\ref{lem:TailProb}.
\end{proofoflemma} 

We proceed to upper bound the summation of the probabilities of the events $\left\{\mu_{i^{*}}<\mu_{i}+2C_{i}^{k}(t)\right\}$ for $t\in\{1,2,\ldots,T\}$ as follows. Using the equation (\ref{eq:Uncertainity}) we have that the inequality
$
\mu_{i^{*}}<\mu_{i}+2C_{i}^{k}(t)
$
implies
\begin{align*}
% \frac{2(\xi+1)\left(N_{i}^{k}(t)+f(t)\right)}{N_{i}^{k}(t)}\frac{\log t}{N_{i}^{k}(t)}&>\frac{\Delta_i^2}{4\sigma_{i}^2}\\
\frac{\Delta_i^2}{4\sigma_{i}^2}\left(N_i^k(t)\right)^2-2(\xi+1)\log t \left(N_i^k(t)+f(t)\right)&<0.
\end{align*}
This inequality does not hold for $N_{i}^{k}(t)> \eta_i(t)$, where
\begin{align*}
\eta_i(t) =\frac{4\sigma_{i}^{2}(\xi+1)}{\Delta_{i}^{2}} \left(1+\sqrt{1+\frac{\Delta_{i}^{2}}{2\sigma_{i}^{2}(\xi+1)}\frac{f(t)}{\log t}}\right)\log t
.
%N_{i}^{k}(t)>\frac{4\sigma_{i}^{2}(\xi+1)}{\Delta_{i}^{2}} %\left(1+\sqrt{1+\frac{\Delta_{i}^{2}}{2\sigma_{i}^{2}(\xi+1)}\frac{f(t)}{\log t}}\right)\log t.
\end{align*}
Thus we have
\begin{align}
\sum_{t=1}^{T}\mathbb{P}\left(Q_{i}^{k}(t)\geq Q_{i^{*}}^{k}(t),N_{i}^{k}(t)<\eta_i(t)\right)\leq \eta_i(T). \label{eq:eta}
\end{align}
%where $\eta(t)=\frac{4\sigma_{i}^{2}(\xi+1)}{\Delta_{i}^{2}} \left[1+\sqrt{1+\frac{\Delta_{i}^{2}}{2\sigma_{i}^{2}(\xi+1)}\frac{f(t)}{\log t}}\right]\log t
%.$\\

Using probability bounds given in Lemma~\ref{lem:TailProb} and equation (\ref{eq:eta}) we prove our main result on performance bounds. 
\begin{theorem}\label{thm:Regret}
The total expected number of times agent $k$ samples suboptimal option $i$ until  time $T$ is upper bounded as
\begin{align*}
&\mathbb{E}\left(n_{i}^{k}(T)\right)\leq   \Gamma(\zeta,\xi,K)+\frac{1}{T^{\xi-1}\log \zeta}\left(\frac{\log K}{T\xi}+\frac{1}{\xi-1}\right)\\
&+  \frac{4\sigma_{i}^{2}(\xi+1)}{\Delta_{i}^{2}} \left(1+\sqrt{1+\frac{\Delta_{i}^{2}}{2\sigma_{i}^{2}(\xi+1)}\frac{f(T)}{\log T}}\right)\log T
\end{align*}
where  $\Gamma(\zeta,\xi,K)=\frac{1}{\log \zeta}(1+\log K)+\frac{1}{2^{\xi}\log \zeta}\left(\frac{\log K}{\xi}+\frac{2}{\xi-1}\right)$ and $\zeta,\xi>1.$
\end{theorem}

{\em Proof of Theorem~\ref{thm:Regret}}: See Appendix.
%The proof of Theorem~\ref{thm:Regret} is given in the appendix. 

Since $f(T)$ is a sublogarithmic function, the expected number of suboptimal samples are logarithmically bounded.

%This concludes the proof of Theorem~\ref{thm:Regret}

\subsection{Performance Measure}\label{SubSecn:Perfom}
In this section we provide a measure to rank the agents according to their relative performance. We motivate with the following two cases for a set of four agents where there is an underlying all-to-all observation structure: 
\begin{center}
\begin{tabular}{*{6}{|c}|}
\hline
& $k$ & 1 & 2 & 3 & 4\\
\hline
Case 1 & $p_{k}$ & 0.5 & 0 & 0 & 0\\
\hline
Case 2 & $p_{k}$ & 0.5 & 1 & 1 & 1\\
\hline
\end{tabular}.
\end{center} 

\begin{figure}[h]
    \centering
    \includegraphics[width=0.4\textwidth]{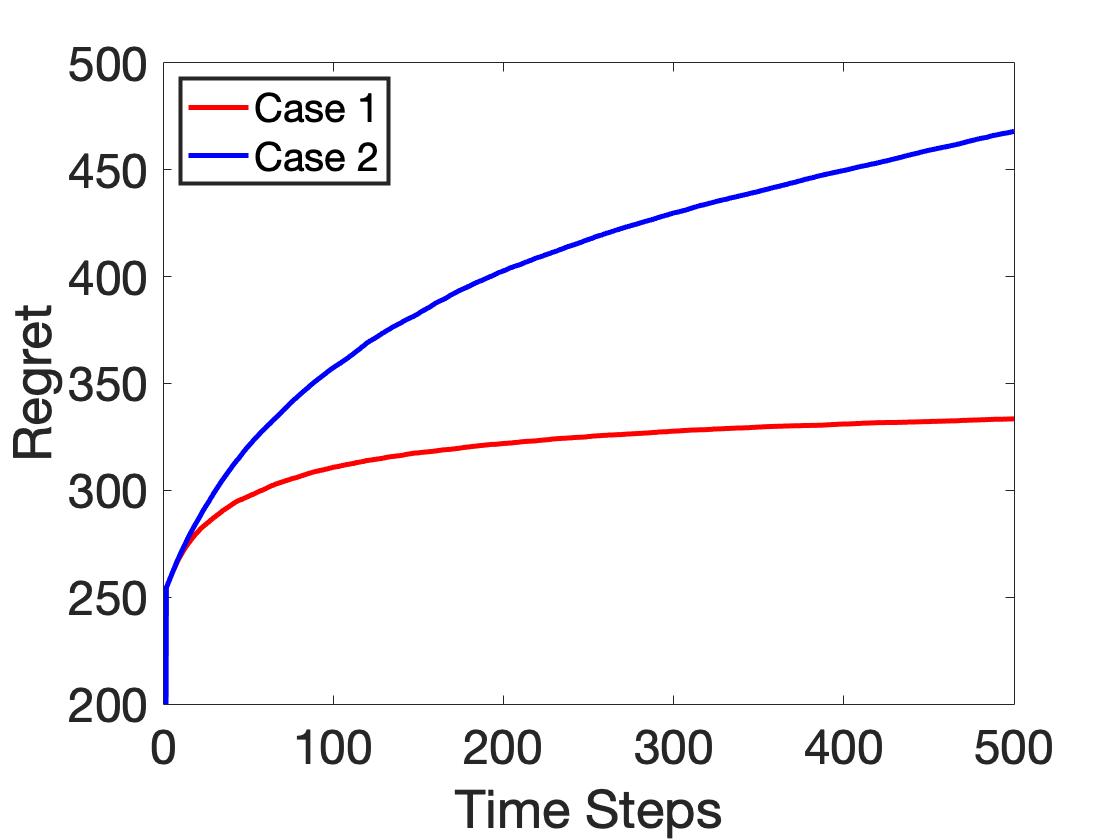}
    \caption{Expected cumulative regret of agent 1 in Cases 1 and 2.}
    \label{Fig:4agent}
\end{figure}

In Case 1, neighbors of agent 1 are not at all sociable and in Case 2 they are maximally sociable. Therefore neighbors of agent 1 explore more in Case 1 than in Case 2. As a result, in Case 1, agent 1 tends to obtain observations from neighbors about lesser known options and this allows agent 1 to exploit more. In Case 2, agent 1 tends to obtain observations from neighbors about well sampled options, and this forces agent 1 to explore more. As a result agent 1 performs better in Case 1 as illustrated in Figure \ref{Fig:4agent}.  

%\begin{multicols}{1}
%\begin{center}
%\begin{figure*}[ht]
%   \includegraphics[width=0.45\textwidth]{4agent.jpg}]
%   \captionof{figure}{Expected cumulative regret of the agent 1}
%   \label{Fig:4agent} 
% \end{figure*}
%\end{center}
%\end{multicols}

With this intuition we propose a measure as follows. First we restrict our attention to a class of problems where the underlying observation structure of the agents is a symmetric $d-$regular graph. This means that $d_{k}=d,\:\forall k$, i.e., every agent has the same number of neighbors. The  all-to-all graph is a special case of the class of $d-$regular graphs with $d=K-1.$
%We also focus on the case in which sociability  $p_k$ takes a value in the discrete set:
We also assume
%\begin{align*}
$p_{k}\in\{.05, .10, .15,\ldots ,1\},\forall k.$
%\end{align*}

We define  performance measure $\epsilon_{p}^{k}\in (0,1]$ for  agent $k$ as
\begin{align}
\epsilon_{p}^{k}&\triangleq \frac{1}{p_{k}+1}\sqrt{\frac{1}{d_{k}}{\sum_{j=1,{j\in\mathcal{V}_{d_{k}}}}^{K}p_{j}}} \; .\label{eq:PM}
\end{align}
%Note that $0<\epsilon_{p}^{k}\leq 1,\:\: \forall k.$ 
Our goal is to show that a lower $\epsilon_p^k$ implies a lower cumulative regret and therefore higher performance % in terms of accumulated reward 
for agent $k$. 
The measure $\epsilon_p^k$ is inversely related to agent $k$'s sociability $p_k$ and directly related to the sociability $p_j$ of the neighbors $j$ of agent $k$. It then makes sense intuitively that lower $\epsilon_p^k$ implies higher performance for agent $k$, since the higher the $p_k$ the more agent $k$ observes and the lower the $p_j$ the more its neighbors explore and the more valuable is their information. 

We next design a protocol for each agent $k$ that depends on $\epsilon_p^k$ and show in Corollary~\ref{cor-regret} that the bound on agent $k$'s cumulative regret is directly related to $\epsilon_p^k$, which suggests that the ordering of agents by $\epsilon_p^k$ predicts the ordering of agents by performance.  We plan to prove in a future publication that the same ordering is predicted even when the protocol does not depend on $\epsilon_p^k$. Simulations in Section~\ref{Secn:Simu} provide validation of these assertions.
%We provide some validation of these assertions with simulations in the next section that these predictions 

%Note that the all-to-all graph is a special case of the class of $d-$regular graphs with $d=K-1.$
Let $f(t)=\epsilon_{p}^{k}.$ Then the objective function becomes
\begin{align*}
	Q^{k}_i(t)&\triangleq \hat{\mu}_{i}^{k}(t)+\sigma_{i}\sqrt{\frac{2(\xi+1)\left(N_{i}^{k}(t)+\epsilon_{p}^{k}\right)}{N_{i}^{k}(t)}\frac{\log t}{N_{i}^{k}(t)}}.\label{eq:UCBQ}
\end{align*}
%For an agent to use this sampling rule, it needs to know %the local probability structure, specifically 
This assumes that each agent $k$ knows the $p_j$ of each of its neighbors $j$.
From Theorem~\ref{thm:Regret}, the
expected number of times agent $k$ samples the suboptimal option $i$ is bounded as
\begin{align*}
&\mathbb{E}\left(n_{i}^{k}(T)\right) \leq  \Gamma(\zeta,\xi,K)
+\frac{1}{T^{\xi-1}\log \zeta}\left(\frac{\log K}{T\xi}+\frac{1}{\xi-1}\right)\\
&+  \frac{4\sigma_{i}^{2}(\xi+1)}{\Delta_{i}^{2}} \left(1+\sqrt{1+\frac{\Delta_{i}^{2}}{2\sigma_{i}^{2}(\xi+1)}\frac{\epsilon_{p}^{k}}{\log T}}\right) \!\log T.
\end{align*}
This suggests that lower $\epsilon_{p}^{k}$ values correspond to lower cumulative expected regret and hence better performance.
 Using the bound on $\mathbb{E}\left(n_{i}^{k}(T)\right)$ and equation (\ref{eq:regretValue}) we upper bound the cumulative expected regret of agent $k$ as follows.
\begin{corollary}
\label{cor-regret}
Let $R^{k}(T)$ be the regret of  agent $k$ up to  time $T$. Then we have
\begin{align*}
&\mathbb{E}\left(R^{k}(T)\right)\leq  \sum_{i=1}^{N}\Delta_{i}\Gamma(\zeta,\xi,K)\\
&+\sum_{i=1}^{N}\left( \frac{\Delta_{i}}{T^{\xi-1}\log \zeta}\left(\frac{\log K}{T\xi}+\frac{1}{\xi-1}\right)\right)\\
&+\sum_{i=1}^{N}\frac{4\sigma_{i}^{2}(\xi+1)}{\Delta_{i}} \left(1+\sqrt{1+\frac{\Delta_{i}^{2}}{2\sigma_{i}^{2}(\xi+1)}\frac{\epsilon_{p}^{k}}{\log T}}\right)\log T.
\end{align*}
Expected cumulative regret is then logarithmically bounded.
\end{corollary}

 The performance measure $\epsilon_p^k$ can be bounded as
\begin{align*}
\epsilon_{p}^{k}\leq \frac{1}{p_{k}+1},\:\: \forall k.
\end{align*}
Accounting for the sociability of neighbors, as we have done, provides a more accurate and tighter bound than only considering individual sociability values. However, for an all-to-all observation structure the rank order can be predicted using only individual sociability values as shown next.% in the following lemma.
\begin{lemma}
Let $\mathcal{G}$  be an all-to-all graph. Let $p_{k}$, $p_{k^{\prime}}$ be the sociability of agents $k$, $k^{\prime}$ such that $p_{k}>p_{k^{\prime}}.$ Then $\epsilon_{p}^{k}<\epsilon_{p}^{k^{\prime}}.$
\end{lemma}
\begin{proofoflemma}
Since $p_{k}>p_{k^{\prime}}$ and $d_k=d_{k^{\prime}}=K-1$ we have 
\begin{align*}
\frac{1}{p_{k}+1}<\frac{1}{p_{k^{\prime}}+1};\:\: \sqrt{\frac{1}{d_k}{\sum_{j=1,{j\in\mathcal{V}_{d_{k}}}}^{K} \!\!p_{j}}}<\sqrt{\frac{1}{d_{k^{\prime}}}{\sum_{j=1,{j\in\mathcal{V}_{d_{k^{\prime}}}}}^{K}\!\! p_{j}}}
\end{align*}
By equation (\ref{eq:PM}), this proves that $\epsilon_{p}^{k}<\epsilon_{p}^{k^{\prime}}.$
\end{proofoflemma}

%%%%%%%%%%%%%%%%%%%%%%%%%%%%%%%%%%%%%%%
\section{Numerical Simulations}\label{Secn:Simu}
We ran numerical simulations to evaluate the performance of  sampling rule   %equations 
(\ref{eq:UCBallocation})--(\ref{eq:UCBQ}) with $f(t) = \epsilon_p^k$ for agent $k$ and with $f(t) = \log\log t$. The results in both cases verify the accuracy of agent ranks predicted by the performance measure $\epsilon_p^k$: lower $\epsilon_p^k$ corresponds to lower cumulative regret, hence higher performance. We show plots in the case $f(t) = \epsilon_p^k$.

We consider 6 agents playing 10-armed bandit problems with two distinct observation structures: A) an all-to-all graph and B) a cyclic graph. In all simulations we let the reward distributions be Gaussian with variance $(\sigma_i^{\prime})^2 = 25$, $i=1, \ldots, 10$, mean values given by
\begin{center}
\begin{tabular}{*{11}{|c}|}
\hline
$i$ & 1 & 2 & 3 & 4 & 5 & 6 & 7 & 8 & 9 & 10\\
\hline
$\mu_i$ & 40 & 50 & 50 & 60 & 70 & 70 & 80 & 90 & 92 & 95 \\
\hline
\end{tabular},
\end{center} 
and sociability values given by
\begin{center}
\begin{tabular}{*{7}{|c}|}
\hline
$k$ & 1 & 2 & 3 & 4 & 5 & 6\\
\hline
$p_{k}$ & 0.50 &0.85 &0.05 &0.50 & 1.00 & 0.90\\
%\cline{2-6}
\hline 
\end{tabular}.
\end{center} 
We provide results for 500 time steps with 1000 Monte Carlo simulations. We set the sampling rule parameter $\xi= 1.1$.

\subsection{All-to-all observation}\label{SubSecn:AlltoAll}
The underlying observation graph structure is all-to-all, equivalently a 5-regular graph. %Since the number of agents is 6, this is essentially an all-to-all graph. 
We calculate the performance measure $\epsilon_p^k$ for each agent $k$ using equation (\ref{eq:PM}):
\begin{center}
\begin{tabular}{*{7}{|c}|}
\hline
$k$ & 1 & 2 & 3 & 4 & 5 & 6\\
\hline
$\epsilon_{p}^{k}$ &0.542 &   0.415  &  0.825 &  0.542 &   0.374  &  0.401\\
\hline 

\end{tabular}.
\end{center} 

\begin{figure}[t]
    \centering
    \includegraphics[width=0.4\textwidth]{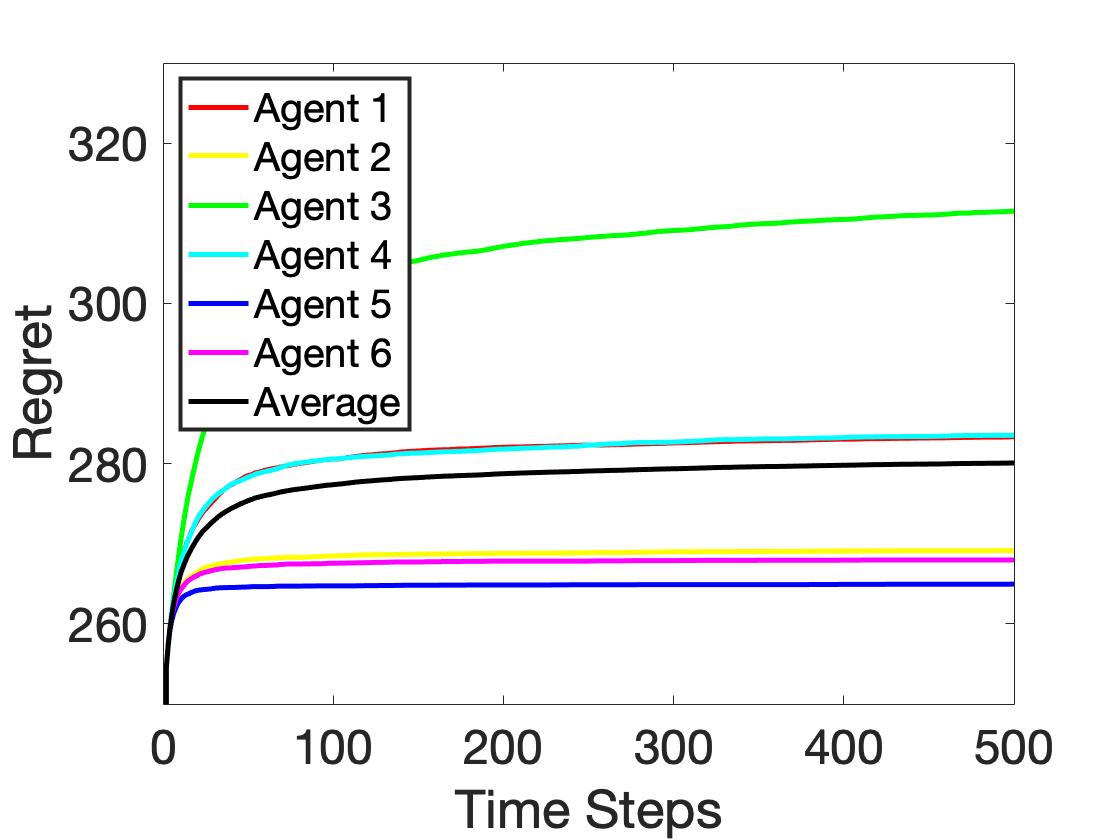}
    \caption{Expected cumulative regret of the 6 agents using the sampling rule given in (\ref{eq:UCBallocation})--(\ref{eq:UCBQ}) and performance measure defined in (\ref{eq:PM}) with distinct observation probabilities $p_k$ and underlying all-to-all observation structure.}
\label{Fig:AlG}
\end{figure}

% \begin{figure*}[!t]
%  \centering
%  \begin{tabular}{ccc}
%   \includegraphics[width=0.45\textwidth]{MeasureAllR.jpg}&   \includegraphics[width=0.45\textwidth]{MeasureR.jpg}\\
%  (a)  Results for all-to-all graph &(b) Results for Cyclic graph
%     \end{tabular}
%   \caption{Expected cumulative regret of the agents, who utilize the sampling rule given in equation (\ref{eq:UCBallocation})--(\ref{eq:UCBQ}) and performance measure defined in equation (\ref{eq:PM}) with distinct observation probabilities.}
%   \label{Fig:AlG} 
% \end{figure*}

The best predicted performer is agent 5 with lowest performance measure  $\epsilon_{p}^{5}=0.374$ and highest sociability value $p_5=1.0$. Second and third best predicted performers are agent 6 and agent 2, respectively, with second and third highest sociability. Agents 1 and 4 are ranked next with equal sociability and $\epsilon_{p}^{1}=\epsilon_{p}^{4}.$ The worst predicted performer is agent 3, with lowest sociability and highest performance measure $\epsilon_{p}^{3}$. %This corresponds to the agent with lowest sociability value and he has the highest cumulative regret. 
These predictions on performance ranking are verified in the simulation results of Figure \ref{Fig:AlG}. 
% Figures \ref{Fig:AlG}(a), \ref{Fig:AlG}(b) and \ref{Fig:AlG}(c) present simulations results for $f(t)=\epsilon_{p}^{k},f(t)=\log\log t$, and $f(t)=0$, respectively. 
The results also verify that for an underlying 
%that when the underlying observation structure is 
all-to-all graph,
the performance rank ordering is predicted by the sociability rank ordering.

\subsection{Cyclic observation}\label{SubSecn:Cyclic}
The underlying observation graph structure is a cycle, equivalently a 2-regular graph, defined %for $k=2,\ldots,5$ agents 
as $ \{k,j\}\in\mathcal{E}$, where $|(k-j)\mod 6|=1$.  We calculate the $\epsilon_p^k$ using (\ref{eq:PM}):
\begin{center}
\begin{tabular}{*{7}{|c}|}
\hline
$k$ & 1 & 2 & 3 & 4 & 5 & 6\\
\hline 
$\epsilon_{p}^{k}$ &0.624  &  0.284  &  0.783  &  0.483  &  0.418  &  0.456\\
\hline 
\end{tabular}.
\end{center}

\begin{figure}[t]
    \centering
    \includegraphics[width=0.4\textwidth]{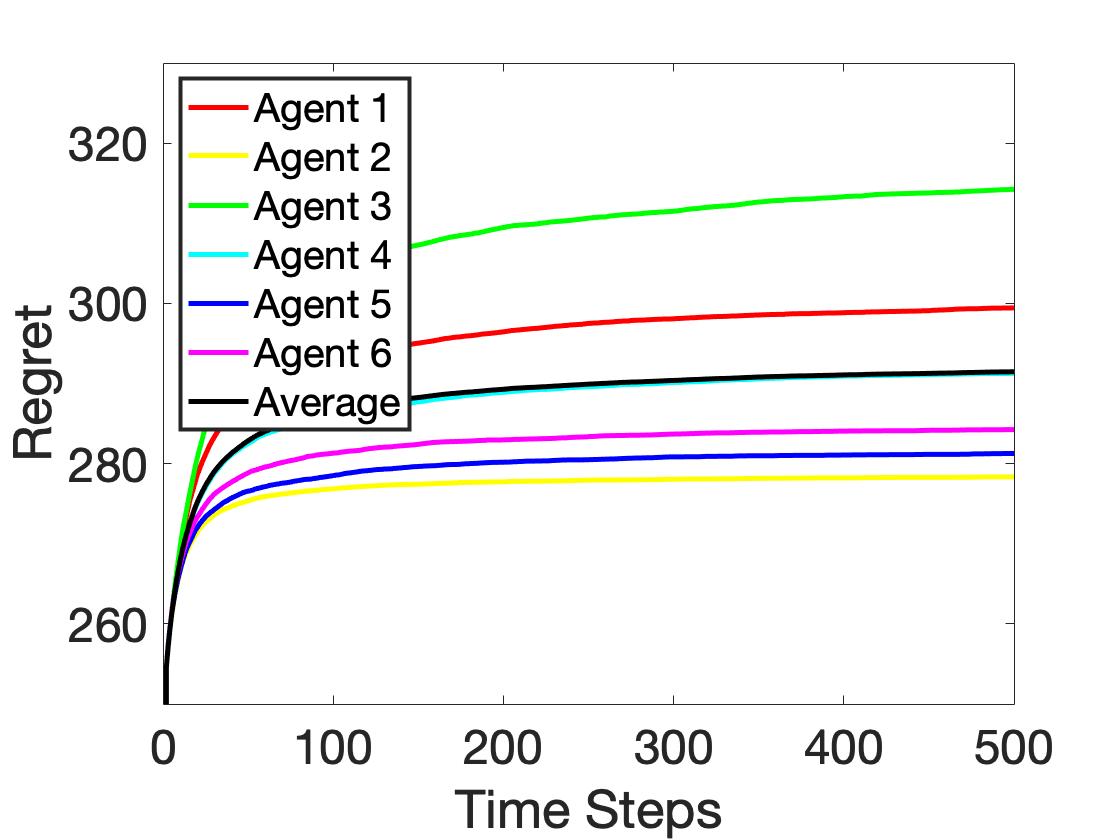}
    \caption{Expected cumulative regret of the 6 agents using the sampling rule given in (\ref{eq:UCBallocation})--(\ref{eq:UCBQ}) and performance measure defined in (\ref{eq:PM}) with distinct observation probabilities $p_k$ and underlying cyclic observation structure.}
\label{Fig:CycG}
\end{figure}

The best predicted performer is agent 2 with the lowest  performance measure  $\epsilon_{p}^{2}=0.284,$ but not the highest sociability. %Recall that sociability values of 
In fact, agents 5 and 6 have higher sociability than agent 2. However, while all three agents 2, 5, and 6 have one neighbor with sociability 0.5, the other neighbor of agents 5 and 6 has sociability 0.9 and 1, respectively, whereas the other neighbor of agent 2 has sociability 0.05. The very low sociability of one of agent 2's neighbors improves agent 2's performance significantly enough that it  outperforms agents 5 and 6. This result illustrates the important role of sociability of neighbors in an agent's performance. Further, while agents 1 and 4 are indistinguishable by their own sociability,
%and they are indistinguishable by considering individual sociability values only. However 
agent 4 has a neighbor with sociability 0.05 whereas agent 1 has neighbors with relatively high sociability. Therefore as predicted by the performance measure, agent 4 outperforms agent 1. Figure \ref{Fig:CycG} validates the predicted rankings.

%%%%%%%%%%%%%%%%%%%%%%%%%%%%%%%%%%%%%%%

\section{Conclusions}\label{Secn:Concl}
We studied a MAMAB problem where agents observe instantaneous actions and rewards of their neighbors according to a stochastic network graph in which agents are distinguished by their sociability, defined as the probability of observing their neighbors. We derived an upper bound for expected cumulative regret of agents.  We proposed a measure to predict relative performance ranking of the agents as a function of sociability of agents and their neighbors. We verified that having less sociable neighbors improves the performance of agents. Accuracy of the measure has been verified analytically through expected cumulative regret bounds and computationally through numerical simulations.

%%%%%%%%%%%%%%%%%%%%%%%%%%%%%%%%%%%%%%%

\section*{acknowledgment}
The first author wishes to thank Peter Landgren for helpful comments during the preparation of this paper.
%%%%%%%%%%%%%%%%%%%%%%%%%%%%%%%%%%%%%%%

\begin{appendix}
\setcounter{proofoftheorem}{0}
\begin{proofoftheorem}
Since $X_i$ is a sub-Gaussian random variable with variance proxy $(\sigma_i^{\prime})^2,$ we have
\begin{align*}
E\left(\exp(\lambda(X_i-\mu_i))\right)\leq \exp\left(\frac{\lambda^2 (\sigma_{i}^{\prime})^2}{2}\right).
\end{align*}

Define a new random variable such that
\begin{align*}
Y_i^k(\tau)&=\sum_{j=1}^{K}\left(X_i^\tau-\mu_i\right )\mathbb{I}_{\{\varphi_j^{\tau} =i\}}\mathbb{I}^{\tau}_{\{k,j\}}.
\end{align*}
Note that $E(Y_i^k(\tau))=E(Y_i^k(\tau)|\mathcal{F}_{\tau-1})=0.$ Let $Z_{i}^k(t)=\sum_{\tau=1}^{t}Y_i^k(\tau).$ 
Let $T_i^k(\tau)=\sum_{j=1}^{K}\mathbb{I}_{\{\varphi_j^{\tau} =i\}}\mathbb{I}^{\tau}_{\{k,j\}}$ For any $\lambda>0$
\begin{align*}
E\left(\exp(\lambda Y_i^k(\tau))|\mathcal{F}_{\tau-1}\right)\leq  \exp\left(\frac{\lambda^2 \sigma_{i}^2}{2}T_i^k(\tau)\right).
\end{align*}
%Recall that 
$\mathbb{I}_{\{\varphi_j^{\tau}=i\}}$ is an $\mathcal{F}_{\tau-1}$ measurable random variable, and so
\begin{align*}
E\left(\exp\left(\lambda Y_i^k(\tau)-\frac{\lambda^2 \sigma_{i}^2}{2}T_i^k(\tau)\right)\Big |\mathcal{F}_{\tau-1}\right)\leq 1.
\end{align*}

Further, using the properties of conditional expectations 
% \begin{align*}
% & E\left(\exp\left(\sum_{1\leq \tau\leq t}\lambda Y_i^k(\tau)-\frac{\lambda^2 \sigma_{i}^2}{2}\sum_{j=1}^{K}\mathbb{I}_{\{\varphi_j^{\tau} =i\}}\mathbb{I}^{\tau}_{\{k,j\}}\right)\right)\\
% &=E\left(E\left(\exp\left(\sum_{1\leq \tau\leq t}\lambda Y_i^k(\tau)-\frac{\lambda^2\sigma_{i}^2K}{2}\sum_{j=1}^{K}\mathbb{I}_{\{\varphi_j^{\tau} =i\}}\mathbb{I}^{\tau}_{\{k,j\}}\right)\Big|\mathcal{F}_{t-1}\right)\right).
% \end{align*}
% Note that
\begin{align*}
& E\left(\exp\left(\sum_{1\leq \tau\leq t}\lambda  Y_i^k(\tau)-\frac{\lambda^2\sigma_{i}^2}{2}T_i^k(\tau)\right)\Big|\mathcal{F}_{t-1}\right)\\
&\leq \exp\left(\sum_{1\leq \tau\leq t-1}\lambda  Y_i^k(\tau)-\frac{\lambda^2\sigma_{i}^2}{2}T_i^k(\tau)\right).
\end{align*}
Thus we see that
\begin{align*}
E\left(\exp\left(\lambda Z_i^k(t)-\frac{\lambda^2 \sigma_{i}^2}{2}N_i^{k}(t)\right)\right) \leq 1.
\end{align*}

The rest of the proof closely follows the papers \cite{Garivier2011,Garivier2008}. For clarity and completeness we include the main steps. % as follows.
Let $\zeta>1.$ Then  $
1\leq N_{i}^{k}(t)\leq \zeta^{D_{t}}$
where $D_{t}=\frac{\log Kt}{\log \zeta}.$ For $\lambda_{j}=\frac{2}{\sigma_{i}}\sqrt{\frac{\kappa\vartheta}{\zeta^{j-1/2}}}$ and $\zeta^{j-1}\leq N_{i}^{k}(t)\leq \zeta^{j}$ we have
\begin{align*}
\frac{2\kappa\vartheta}{\lambda_{j}\sqrt{N_{i}^{k}(t)}}+\frac{\sigma_{i}^{2}}{2}\lambda_{j}\sqrt{N_{i}^{k}(t)}\leq \sqrt{\vartheta},
\end{align*}
where $\kappa=\frac{1}{\sigma_i^2\left(\zeta^{\frac{1}{4}}+\zeta^{-\frac{1}{4}}\right)^2}.$

Recall  from the Markov inequality that 
%\begin{align*}
$\mathbb{P}(Y \geq a )\leq \frac{\mathbb{E}(Y)}{a}$
%\end{align*}
for any random variable $Y$. Thus,
\begin{align*}
\mathbb{P}\left(\frac{Z_{i}^{k}(t)}{\sqrt{N_{i}^{k}(t)}}\geq\sqrt{ \vartheta}\right)\leq \cup_{j=1}^{D_{T}}\exp(-2\kappa\vartheta).
\end{align*}
\end{proofoftheorem}
%This concludes the proof of Theorem \ref{thm:TailProb}.
%%%%%%%%%%%%%%%%%%%%%%%%%%
\setcounter{proofoftheorem}{1}
\begin{proofoftheorem}
From equations (\ref{eq:Prob})--(\ref{eq:eta}) we have
\begin{align*}
&\mathbb{E}\left(n_{i}^{k}(T)\right)
\leq \sum_{t=1}^{T}\mathbb{P}\left(\hat{\mu}_{i^{*}}^{k}\leq \mu_{i^{*}}-C_{i^{*}}^{k}(t),N_{i}^{k}(t)>\eta_i(t)\right)\\
&+\eta_i(T)+\sum_{t=1}^{T}\mathbb{P}\left(\hat{\mu}_{i}^{k}\geq \mu_{i}+C_{i}^{k}(t),N_{i}^{k}(t)>\eta_i(t)\right).
\end{align*}
Note that $\frac{N_{i}^{k}(t)+f(t)}{N_{i}^{k}(t)}\geq 1,\forall t\geq 1$. Using Lemma~\ref{lem:TailProb} and  $\delta(\xi)=\xi+1-\delta^{\prime}(\epsilon)$  we have
\begin{align*}
\mathbb{E}\left(n_{i}^{k}(T)\right)
\leq \eta_i(T) +2\nu\sum_{t=1}^{T}\frac{\log Kt}{t^{\delta(\xi)}}.
\end{align*}
Since $\delta^{\prime}(\epsilon)$ can be made arbitrarily small, the summation term can be evaluated as follows:
\begin{align*}
\sum_{t=1}^{T}\frac{\log Kt}{t^{\delta(\xi)}}\leq& \sum_{t=1}^{T}\frac{\log K}{t^{\xi+1}}+\sum_{t=1}^{T}\frac{1}{t^{\xi}}\\
\leq& \log K\left(1+\int_{2}^{T}\frac{1}{t^{\xi+1}}dt\right)+1+\int_{2}^{T}\frac{1}{t^{\xi}}dt\\
=& \; 1+\log K+\frac{1}{2^{\xi}}\left(\frac{\log K}{\xi}+\frac{2}{\xi-1}\right)\\
&+\frac{1}{T^{\xi}}\left(\frac{\log K}{\xi}+\frac{T}{\xi-1}\right).
\end{align*}
Thus we have 
\begin{align*}
&\mathbb{E}\left(n_{i}^{k}(T)\right)\leq  \frac{1}{\log \zeta}(1+\log K)+\frac{1}{2^{\xi}\log \zeta}\left(\frac{\log K}{\xi}+\frac{2}{\xi-1}\right)\\
&\:\:\:+\frac{1}{T^{\xi-1}\log \zeta}\left(\frac{\log K}{T\xi}+\frac{1}{\xi-1}\right)\\
&\:\:\:+\frac{4\sigma_{i}^{2}(\xi+1)}{\Delta_{i}^{2}} \left(1+\sqrt{1+\frac{\Delta_{i}^{2}}{2\sigma_{i}^{2}(\xi+1)}\frac{f(T)}{\log T}}\right)\log T.
\end{align*}
%which concludes the proof of Theorem~\ref{thm:Regret}.
\end{proofoftheorem}
\end{appendix}

%%%%%%%%%%%%%%%%%%%%%%%%%%

%    \input{MAMAB.bbl}
\bibliographystyle{IEEEtran}
\bibliography{MAMAB}

\begin{thebibliography}{10}
\providecommand{\url}[1]{#1}
\csname url@samestyle\endcsname
\providecommand{\newblock}{\relax}
\providecommand{\bibinfo}[2]{#2}
\providecommand{\BIBentrySTDinterwordspacing}{\spaceskip=0pt\relax}
\providecommand{\BIBentryALTinterwordstretchfactor}{4}
\providecommand{\BIBentryALTinterwordspacing}{\spaceskip=\fontdimen2\font plus
\BIBentryALTinterwordstretchfactor\fontdimen3\font minus
  \fontdimen4\font\relax}
\providecommand{\BIBforeignlanguage}[2]{{%
\expandafter\ifx\csname l@#1\endcsname\relax
\typeout{** WARNING: IEEEtran.bst: No hyphenation pattern has been}%
\typeout{** loaded for the language `#1'. Using the pattern for}%
\typeout{** the default language instead.}%
\else
\language=\csname l@#1\endcsname
\fi
#2}}
\providecommand{\BIBdecl}{\relax}
\BIBdecl

\bibitem{Sutton}
R.~S. Sutton and A.~G. Barto, \emph{Introduction to Reinforcement
  Learning}.\hskip 1em plus 0.5em minus 0.4em\relax MIT Press Cambridge, MA,
  USA, 1998.

\bibitem{Robbins}
H.~Robbins, \emph{Some Aspects of the Sequential Design of Experiments}.\hskip
  1em plus 0.5em minus 0.4em\relax Springer New York, 1985.

\bibitem{Gittins}
J.~C. Gittins, ``Bandit processes and dynamic allocation indices,''
  \emph{Journal of the Royal Statistical Society. Series B (Methodological)},
  vol.~41, pp. 148--177, 1979.

\bibitem{LaiRobbins}
T.~L. Lai and H.~Robbins, ``Asymptotically efficient adaptive allocation
  rules,'' \emph{Advances in Applied Mathematics}, vol.~6, no.~1, pp. 4--22,
  1985.

\bibitem{AgrawalSimpl}
R.~Agrawal, ``Sample mean based index policies with o(log n) regret for the
  multi-armed bandit problem.'' \emph{Advances in Applied Probability},
  vol.~27, pp. 1054--1078, 1995.

\bibitem{Auer}
P.~Auer, N.~Cesa-Bianchi, and P.~Fisher, ``Finite-time analysis of the
  multi-armed bandit problem.'' \emph{Machine Learning}, vol.~47, pp. 235--256,
  2002.

\bibitem{Kauffman}
E.~Kauffman, O.~Cappe, and A.~Garivier, ``On {B}ayesian upper confidence bounds
  for bandit problem,'' in \emph{International Conference on Artificial
  Intelligence and Statistics,}, April 2012, pp. 592--600.

\bibitem{Reverdy}
P.~Reverdy, V.~Srivastava, and N.~E. Leonard, ``Modeling human decision-making
  in generalized {Gaussian} multi-armed bandits,'' in \emph{Proceedings of the
  IEEE}, vol. 102, no.~4, 2014, pp. 544--571.

\bibitem{1104491}
V.~Anantharam, P.~Varaiya, and J.~Walrand, ``Asymptotically efficient
  allocation rules for the multiarmed bandit problem with multiple plays-part
  i: I.i.d. rewards,'' \emph{IEEE Transactions on Automatic Control}, vol.~32,
  no.~11, pp. 968--976, 1987.

\bibitem{PeterCDC}
P.~Landgren, V.~Srivastava, and N.~E. Leonard, ``Distributed cooperative
  decision-making in multiarmed bandits: Frequentist and {Bayesian}
  algorithms,'' in \emph{IEEE Conference on Decision and Control}, December
  2016, pp. 167--172.

\bibitem{6763073}
D.~Kalathil, N.~Nayyar, and R.~Jain, ``Decentralized learning for multiplayer
  multiarmed bandits,'' \emph{IEEE Transactions on Information Theory},
  vol.~60, no.~4, pp. 2331--2345, 2014.

\bibitem{PeterECC}
P.~Landgren, V.~Srivastava, and N.~E. Leonard, ``On distributed cooperative
  decision-making in multiarmed bandits,'' in \emph{European Control
  Conference}, June 2016, pp. 243--248.

\bibitem{Gopalan16}
R.~K. Kolla, K.~Jagannathan, and A.~Gopalan, ``Collaborative learning of
  stochastic bandits over a social network,'' \emph{arXiv:1602.08886v2}, 2016.

\bibitem{LandgrenCDC2018}
P.~Landgren, V.~Srivastava, and N.~E. Leonard, ``Social imitation in
  cooperative multiarmed bandits: partition-based algorithms with strictly
  local information,'' in \emph{IEEE Conference on Decision and Control},
  December 2018, pp. 5239--5244.

\bibitem{Tilman2016MaintainingCI}
A.~R. Tilman, J.~R. Watson, and S.~Levin, ``Maintaining cooperation in
  social-ecological systems:,'' \emph{Theoretical Ecology}, vol.~10, pp.
  155--165, 2016.

\bibitem{Lai}
\BIBentryALTinterwordspacing
T.~L. Lai, ``Adaptive treatment allocation and the multi-armed bandit
  problem,'' \emph{Ann. Statist.}, vol.~15, no.~3, pp. 1091--1114, 09 1987.
  [Online]. Available: \url{https://doi.org/10.1214/aos/1176350495}
\BIBentrySTDinterwordspacing

\bibitem{Garivier2011}
A.~Garivier and E.~Moulines, \emph{On Upper-Confidence Bound Policies for
  Switching Bandit Problems}.\hskip 1em plus 0.5em minus 0.4em\relax Berlin,
  Heidelberg: Springer Berlin Heidelberg, 2011, pp. 174--188.

\bibitem{Garivier2008}
------, ``On upper-confidence bound policies for non-stationary bandit
  problems,'' \emph{arXiv:0805.3415v1}, 2008.

\end{thebibliography}

\newpage
\end{document}